\documentclass[10pt]{amsart}
\usepackage{amsmath,amssymb,amsfonts} % Typical maths resource packages
\usepackage{graphics}                 % Packages to allow inclusion of graphics
\usepackage[square, numbers]{natbib}
\usepackage[english]{babel}
\usepackage{tikz} 
\usepackage[utf8]{inputenc}
\usepackage[all]{xy}	
\usepackage[T1]{fontenc}
\usepackage{pifont}

\usepackage[margin=1.2in]{geometry}

\newcommand{\R}{\mathbb{R}}

\theoremstyle{definition}
\newtheorem{defn}{Definition}[subsection]
\newtheorem{rem}[defn]{Remark}
\newtheorem{cons}[defn]{Construction}
\newtheorem{exam}[defn]{Example}

\newtheorem{nota}[defn]{Notation}

\theoremstyle{plain}
\newtheorem{Prop}[defn]{Proposition}
\newtheorem{Thm}[defn]{Theorem}
\newtheorem{cor}[defn]{Corollary}
\newtheorem{lem}[defn]{Lemma}

\theoremstyle{exampstyle}\newtheorem{thmx}{Theorem}

\theoremstyle{exampstyle}\newtheorem{thmxbis}{Theorem}

\def\Ind{\setbox0=\hbox{$x$}\kern\wd0\hbox to 0pt{\hss$\mid$\hss}
\lower.9\ht0\hbox to 0pt{\hss$\smile$\hss}\kern\wd0}
\def\Notind{\setbox0=\hbox{$x$}\kern\wd0\hbox to 0pt{\mathchardef
\nn=12854\hss$\nn$\kern1.4\wd0\hss}\hbox to
0pt{\hss$\mid$\hss}\lower.9\ht0 \hbox to
0pt{\hss$\smile$\hss}\kern\wd0}

\newcommand{\F}{\mathcal{F}}

\newtheoremstyle{exampstyle}
{3pt} % Space above
{3pt} % Space below
{\itshape} % Body font
{} % Indent amount
{\bfseries} % Theorem head font
{.} % Punctuation after theorem head
{.5em} % Space after theorem head
{} % Theorem head spec (can be left empty, meaning `normal')

\makeindex

\begin{document}

\title[Differential fields and geodesic flows II]{Differential fields and Geodesic flows II :\\  \footnotesize\mdseries Geodesic flows of pseudo-Riemannian algebraic varieties}
\author{Rémi Jaoui}
\address{Rémi Jaoui, Département de Mathématiques, Université Paris-Sud, Bâtiment 425, 91405 Orsay, France}
\email{remi.jaoui@math.u-psud.fr}
\date\today

\maketitle

\begin{abstract}
We define the notion of a smooth pseudo-Riemannian algebraic variety $(X,g)$ over a field $k$ of characteristic $0$, which is an algebraic analogue of the notion of Riemannian manifold and we study, from a model-theoretic perspective, the algebraic differential equation describing the geodesics on $(X,g)$.

When $k$ is the field of real numbers, we prove that if the real points of $X$ are Zariski-dense in $X$ and if the real analytification of $(X,g)$ is a compact Riemannian manifold with negative curvature, then the algebraic differential equation describing the geodesics on $(X,g)$ is absolutely irreducible and its generic type is orthogonal to the constants.
\end{abstract}

\vspace{0,5cm}

1. This article is the second one in a series of two articles devoted  to the construction of large families of algebraic ordinary differential equations whose generic types, in the theory of differentially closed fields, are orthogonal to the constants.

In the first part of this work \cite{moi}, we have translated this model-theoretic property of an algebraic differential equation into a property of the associated $D$-variety $(X,v)$. Here $X$ denotes a (smooth) algebraic variety over the base field $k$ and $v$ a regular vector field on $X$, that defines the algebraic differential equation under study, and the orthogonality to the constants becomes ``the $D$-variety $(X,v)^n$ admits \textit{no non-constant rational integral for any natural number} $n$''. Then, when the base field is a subfield of the field $\R$ of real numbers, we have proved that $(X,v)$ satisfies this property as soon as the flow of the dynamical system $(X(\R), v_\R)$ defined by the real points of $(X,v)$ is ``\textit{sufficiently topologically mixing}''.

In this second part, we provide an algebraic framework, complementary to the model-theoretic developments in \cite{moi}, to study geodesic flows of Riemannian manifolds from the perspective of differential algebra and model theory.

In order to associate a definable set in a differentially closed field to a differential equation, we need to express it as \textit{a formula} in the language $\mathcal L_\delta = \lbrace 0,1,+,-,\times,\delta \rbrace$ of differential rings. Consequently, an important limitation for us is to work with \textit{polynomial differential equations}, whereas the setting for classical Riemannian geometry is the category of smooth manifolds (and therefore a broader class of differential equations). Accordingly, we develop a formalism of Riemannian geometry in the framework of smooth algebraic varieties over a field of characteristic $0$.

Once this has been settled, we have a large class -- the geodesic flows of various real algebraic smooth pseudo-Riemannian varieties -- of algebraic differential equations at hand, that we can start to investigate from a model-theoretic perspective. Using the results of the first part of this work, we address the problem of \textit{orthogonality to the constants} for their generic types, under suitable irreducibility assumptions. \\

2. Originally, the \textit{geodesics of a Riemannian manifold $(M,g)$} were defined as the smooth curves drawn on $M$, with constant velocity and minimal length among the curves with the same end points. For a smooth submanifold $M$ of the Euclidean space $\R^n$ (endowed with the restriction $g$ of the Euclidean metric), the geodesics of $(M,g)$ all share the property that \textit{their acceleration is normal to the submanifold $M$.}

In modern terminology, the term of geodesic is understood in a broader meaning, referring to the curves lying on $M$ satisfying the latter property and therefore, to the solutions of a second-order differential equation on $M$. Alternatively, one may define the geodesics on a submanifold $M$ of the Euclidean space $\R^n$, in physical terms, as the trajectories of a particle, constrained to move without friction on the submanifold $M$.

There are two classical ways of writing this differential equation intrinsically, in terms of the Riemannian manifold $(M,g)$. On the one hand, one can use the \textit{Levi-Civita connection} on the manifold $M$ to describe the geodesics as the curves whose covariant acceleration vanishes. It follows from this formulation that the geodesics (lifted canonically to $TM$) define a foliation of $TM$ by curves.

On the other hand, using the basic principles of \textit{Hamiltonian Mechanics}, one can describe directly the vector field $v$ on $TM$ whose integral curves are the lifted geodesics. For instance, the vector field $v$ is the Hamiltonian vector field associated with the free Hamiltonian $H_g(x,v) = \frac 1 2 g_x(v,v)$ (with respect to the symplectic structure on $TM$ induced by the canonical symplectic structure on $T^\ast M$ via the isomorphism $TM \simeq T^\ast M$ defined by $g$). The flow on $TM$ of this vector field is called \textit{the geodesic flow of the Riemannian manifold $(M,g)$}.      

\textit{The integrability properties} of the differential equation describing the geodesics of a Riemannian manifold $(M,g)$ and the related problem of describing the \textit{dynamic of the geodesic flow} have been a main focus for many great mathematicians. Notably, Jacobi explored some highly non-trivial cases of complete integrability for these systems while Hadamard \cite{Had}, following ideas of Poincaré, was the first one to understand that \textit{negative curvature} was the source of non-integrability results for these differential equations. Later, the results of Hadamard were generalised to arbitrary compact Riemannian manifolds with negative curvature by Anosov \cite{Ano} and his followers (see for example \cite{Pan}).

It is natural to expect that, when dealing with algebraic differential equations, the aforementioned results should have important consequences for the model-theoretic behaviour of \textit{the definable set} (in a differentially closed field) \textit{associated to the geodesic differential equation}. 
\\

3. Fix $k$ a field of characteristic $0$. There is a natural analogue -- the notion of \textit{smooth pseudo-Riemannian algebraic variety over $k$} -- for the notion of Riemannian manifold, namely a smooth algebraic variety $X$ over $k$ endowed with a non-degenerate algebraic symmetric $2$-form $g$.

Many of the basic constructions in Riemannian geometry can still be carried out in this context. For instance, we can define the Levi-Civita connection on the sheaf $\Theta_{X/k}$ of vector fields on $X$ and various avatars of the curvature. Moreover, we associate to any smooth pseudo-Riemannian variety over $k$, a $D$-variety $(T_{X/k},v)$ over the differential field $(k,0)$ called \textit{the geodesic $D$-variety of $(X,g)$}.

When $k$ is the field of real numbers, one may consider the real analytification functor which associates to a smooth algebraic pseudo-Riemannian variety $(X,g)$ over $\R$, a (real analytic) pseudo-Riemannian manifold $(X(\R),g_\R)^{an}$. In the second section, we check that our formal definitions and constructions are compatible with the classical ones for the pseudo-Riemannian manifold $(X(\R),g_\R)^{an}$ through real analytification.

In order to gather diverse examples of geodesic $D$-varieties, we need some tools to construct interesting algebraic  pseudo-Riemannian varieties. Working over the field of real numbers, these tools will be approximation theorems asserting that  any compact Riemannian manifold may be realised, after a small perturbation of the metric, as the analytification of some algebraic pseudo-Riemannian variety over $\R$:
{\thmx \label{Int1} Let $(M,g_M)$ be a smooth, connected and compact Riemannian manifold. For every neighbourhood $\mathcal U \subset \mathcal C^\infty(S^2T^\ast M)$ of $g_M$ for the $\mathcal C^\infty$--topology, there exist a smooth pseudo-Riemannian variety $(X,g)$ over $\R$ and a diffeomorphism $\phi : M \longrightarrow X(\R)^{an}$ such that $\phi^\ast g^{an}_{\R} \in \mathcal U$.} 
\vspace{0.2cm}

A consequence of Theorem \ref{Int1} is the following existence principle: for any \textit{property $(P)$} of a compact Riemannian manifold \textit{preserved by small perturbations}, there exists \textit{a smooth pseudo-Riemannian variety over $\R$} whose analytification \textit{satisfies the property $(P)$}. In particular, Theorem \ref{Int1} provides many examples of smooth pseudo-Riemannian varieties $(X,g)$ over $\R$ whose analytification is a compact Riemannian manifold with negative curvature. 

We also prove an ``embedded version'' of Theorem \ref{Int1}, that ensures the existence of many compact and non-singular real algebraic subsets $M$ of the Euclidean space $\R^n$ such that the restriction to $M$ of the Euclidean metric has negative curvature:
  
{\thmxbis\label{Int2} Let $(M,g)$ be a smooth, connected and compact Riemannian manifold. For every neighbourhood $\mathcal U \subset \mathcal C^\infty(S^2T^\ast M)$ of $g_M$ for the $\mathcal C^\infty$--topology, there exist $m \in \mathbb{N}$ and a $\mathcal C^\infty$-embedding $i :  M \longrightarrow \R^m$ such that:
\begin{itemize}
\item[(i)] The image of $M$ is a non-singular algebraic subset $X(\R)$ of $\R^m$, defined by some smooth, locally closed subvariety $X$ of $\mathbb{A}^m_\R$.
\item[(ii)] Let $g_{0}$ be the usual Euclidean metric on $\R^m$. The pullback $i^\ast g_0$ of the Riemannian metric $g_0$ is in $\mathcal U$. 
\end{itemize}}
%A consequence of Theorem \ref{Int2} is that there exist many \textit{compact and non-singular real algebraic subset} $M$ of \textit{the Euclidean space} $\R^n$ such that the restriction to $M$ of the Euclidean metric has \textit{negative curvature}.

The proofs of Theorem \ref{Int1} and of Theorem \ref{Int2} are both based on Tognoli's Theorem, which asserts the existence of algebraic embeddings for smooth compact manifolds. Taking this theorem for granted, Theorem \ref{Int1} is a direct application of the Stone-Weierstrass Theorem, while the proof of Theorem \ref{Int2} is slightly more involved. The key point behind its proof is a clever construction of Nash, known as \textit{Nash twist}, that plays a central role in his existence proofs for isometric embeddings of differentiable Riemannian manifolds (see \cite{Nash2} and \cite{Nash}; see also \cite{Gromov} for a modern presentation). \\ 
 
4. Let $(X,g)$ be an absolutely irreducible pseudo-Riemannian variety over $k$ and $(T_{X/k},v)$ the associated geodesic $D$-variety. The Hamiltonian $H_g$ defining $v$ is \textit{a rational integral} of the vector field $v$ and defines a morphism of $D$-varieties over $(k,0)$: 

$$ H_g : (T_{X/k},v) \longrightarrow (\mathbb{A}^1,0)$$

In particular, no matter which pseudo-Riemannian variety we started from, the generic type of the geodesic $D$-variety is never orthogonal to the constants. Consequently, the interesting property is the \textit{orthogonality to the constants for} (the generic type of) \textit{the fibres of $H_g$}. A consequence of the homogeneity properties of $H_g$ is that this property does not depend on the non zero fibre chosen.

In the case of the sphere $\mathbb{S}^2 \subset \mathbb{R}^3$ endowed with the restriction of the Euclidean metric (seen as a pseudo-Riemannian variety over $\R$), the geodesic $D$-variety admits a complete system of rational integrals and therefore all the fibres of $H_0$ are non-orthogonal to the constants. However, when $k$ is a subfield of $\R$ and when $(X(\R),g_\R)^{an}$ is a compact Riemannian manifold \textit{with negative curvature}, we prove that the opposite situation arises. Indeed, we prove: 

{\thmx\label{Int3} Let $(X,g)$ be an absolutely irreducible smooth pseudo-Riemannian variety over $\R$, with Zariski-dense real points, and such  that $(X(\R),g_\R)^{an}$  is a non empty, compact Riemannian manifold.

If the sectional curvature of this Riemannian manifold is negative, then every non zero fibre of the morphism of $D$-varieties
$$ H_g : (T_{X/k} ,X_g) \longrightarrow (\mathbb{A}^1,0)$$
is an absolutely irreducible $D$-variety with  generic type orthogonal to the constants.
\vspace{0.2cm}}

First note that Theorem \ref{Int3} implies its ``embedded version'' which was stated as Theorem A in \cite{moi} (Corollary \ref{theoremA} in this article).

Theorem \ref{Int3} is proved by a combination of the dynamical criterion of orthogonality to the constants proved in \cite{moi} and of the results of Anosov and his followers on the topological dynamics of the unitary geodesic flow.

In conclusion, let us indicate that we expect a stronger form of Theorem \ref{Int3} to hold. Namely, under the assumption of Theorem \ref{Int3}, \textit{the generic types of the fibres of $H_g$} should be \textit{minimal and trivial types}, at least when $\mathrm{dim}(X) = 2$. \\

5. This article is organized as follows. In the first section, we define the notion of a smooth pseudo-Riemannian variety over a field $k$ of characteristic $0$. We also give analogues in this setting of some basic formal objects of Riemannian geometry. In particular, we associate to any such pseudo-Riemannian variety, a $D$-variety over $(k,0)$ called the geodesic $D$-variety.

In the second section, we study some properties of the analytification functor. We show that the notions of the first section and the classical notions of Riemannian geometry are related through analytification. We also give proofs for Theorem \ref{Int1} and Theorem \ref{Int2}.

In the third section, we give a self-contained exposition of the mixing properties of the geodesic flow of a Riemannian compact manifold with negative curvature, involved in the proof of Theorem \ref{Int3}. There is nothing really original in this section and its main result is well-known by experts of hyperbolic dynamics. However, we need some specific form of these mixing properties, namely  for the topological dynamics of the geodesic flow of a compact Riemannian manifold with negative (but variable) curvature, for which we could not find an explicit statement in the literature, and we explain how to derive it from available references.

In the last section, we first deal with the problem of absolute irreducibility for the fibres of the Hamiltonian integral of the geodesic $D$-variety. Then, we gather the previous constructions and results to complete the proofs of Theorem \ref{Int3} and its corollaries. \\

6. The results of this article constitute a part of my PhD. thesis, supervised by Jean-Benoît Bost (Orsay) and by Martin Hils (Paris VII--Münster) and they are based on their many interesting ideas and suggestions. I would also like to thank Frederic Paulin for discussing with me the dynamical properties of geodesic flows, described in the third part of this article. 

\tableofcontents
\section{Hamiltonian dynamics and pseudo-Riemannian varieties}

We fix $k$ a field of characteristic $0$. By a variety over $k$, we mean a reduced and separated scheme of finite type over $k$.

If $X$ is a variety over $k$, the tangent bundle (resp. cotangent bundle) of $X$ is denoted $T_{X/k}$ (resp. $T_{X/k}^{\ast}$). The coherent sheaf $\Theta_{X/k}$ (resp. $\Omega^1_{X/k}$) of vector fields on $X$ (resp. of $1$-forms on $X$) is the sheaf of sections of $T_{X/k}$ (resp. $T^{\ast}_{X/k}$).

If $X$ is a smooth variety over $k$, then $\Theta_{X/k}$ and $\Omega^1_{X/k}$ are locally free sheaves and the pairing 
\begin{eqnarray}\label{pairing} \Theta_{X/k} \times \Omega^{1}_{X/k} \longrightarrow \mathcal O_X \end{eqnarray} 
gives rise to isomorphisms $\Theta_{X/k} \simeq (\Omega_{X/k}^{1})^\vee$ and $\Omega^1_{X/k} \simeq  \Theta_{X/k}^\vee$.

\subsection{Hamiltonian formalism} 
{\defn A \textit{smooth symplectic variety over $k$} is a pair $(X,\omega)$ where $X$ is a smooth variety over $k$ and $\omega \in \mathrm{H}^0(X,\Lambda^2 \Omega^1_{X/k})$ is a closed non-degenerate alternating $2$-form on $X$.}

{\cons Let $(X,\omega)$ be a smooth symplectic variety over $k$. The $2$-form $\omega$ induces a morphism of coherent sheaves over $X$: 
$$\Phi : \Theta_{X/k} \longrightarrow \Omega^1_{X/k} $$
such that, for any local section $v$ of $\Theta_{X/k}$, $\Phi(v) = \omega(v,.)$.

Because the $2$-form $\omega$ is non-degenerate, this morphism is an isomorphism. 
In particular, we get an isomorphism of $k$-vector spaces
\begin{eqnarray}\label{isomorphism}\Phi : \mathrm{H}^0(X,\Theta_{X/k}) \overset{\sim}{\longrightarrow} \mathrm{H}^0(X,\Omega^1_{X/k}).\end{eqnarray}}

{\defn Let $(X,\omega)$ be a smooth symplectic variety over $k$ and $H \in \mathrm{H}^0(X,\mathcal O_{X/k})$ a function on $X$. The \textit{Hamiltonian vector field $X_H \in \mathrm{H}^0(X,\Theta_{X/k})$ associated to $H$} is the vector field corresponding to the global $1$-form $dH \in \mathrm{H}^0(X,\Omega^1_{X/k})$ in the isomorphism (\ref{isomorphism}).} 

{\lem\label{firstintegral} Let $(X,\omega)$ be a smooth symplectic variety and $H \in \mathrm{H}^0(X,\mathcal O_{X/k})$ a function on $X$. The function $H$ is an integral of the vector field $X_H$.}

By definition, we have $X_H(H) = dH(X_H) = \omega(X_H,X_H) = 0$ because $\omega$  is alternating.

{\exam\label{cotangentbundle} Let $X$ be a smooth variety over $k$. There exists a canonical symplectic structure on the cotangent bundle $T^\ast_{X/k}$ of $X$, that we describe below.

Let $Y = T_{X/k}^\ast$ and $\pi : Y \longrightarrow X$ the canonical projection. First note that, as $X$ is a smooth variety over $k$, the variety $Y$ is also smooth over $k$.

Moreover, the pairing (\ref{pairing}) induces a morphism of coherent sheaves over $X$: 
$$ \tilde{\theta} : \Theta_{X/k} \longrightarrow \pi_\ast \mathcal O_Y \simeq \mathrm{Sym}((\Omega^1_{X/k})^\vee)$$

Using the adjunction between the functors $\pi_\ast$ and $\pi^\ast$ and composing with the differential $d\pi : \Theta_{Y/k} \longrightarrow \pi^\ast \Theta_{X/k}$  of $\pi$, we get a morphism of coherent sheaves over $Y$: 
$$\theta = \tilde{\theta} \circ d \pi : \Theta_{Y/k} \longrightarrow \mathcal O_Y$$

Then, $\theta \in \mathrm{H}^0(Y,\Theta_{Y/k}^\vee) \simeq \mathrm{H}^0(Y,\Omega^1_{Y/k})$ may be identified with a $1$-form on $Y$ called \textit{the canonical $1$-form of $T^\ast_{X/k}$}.

{\rem\label{etale} Let $f : X \longrightarrow X'$ be an étale morphism of smooth varieties over $k$. The morphism $g = T^\ast f : T^\ast_{X/k} \longrightarrow T^\ast_{X'/k}$ is also étale.
Let $\theta_X$ and $\theta_{X'}$ be the canonical $1$-forms of $T^\ast_{X/k}$ and $T^\ast_{X'/k}$. By using the construction above, it is easy to check that $g^\ast\theta_{X'} = \theta_X$.}

Using this remark, we can give the following description of the canonical $1$-form. Let $X$ be a smooth variety over $k$. As $X$ is smooth, there exist étale charts for $X$, namely a covering of $X$ by open subspaces $U \subset X$  and étale morphisms $$f = (x_1, \cdots, x_n) : U  \longrightarrow \mathbb{A}^n.$$

Let $\pi : T^\ast_{X/k} \longrightarrow X$ be the projection. As $f$ is an étale morphism, it defines a trivialisation of the cotangent bundle of $X$:
$$ \pi^{-1}(U) \simeq U \times \mathbb{A}^n$$
whose coordinate functions are denoted $p_1,\ldots, p_n \in \mathcal O_X(U)$.

Using remark \ref{etale} and the direct computation for the canonical $1$-form in $\mathbb{A}^n$ \citep[Chapter 2]{Arn2}, the canonical $1$-form on $X$ is given by: 
\begin{eqnarray}\label{canonicalform} \theta = \sum_{i =1}^n p_i. dx_i.
\end{eqnarray}

{\lem Let $X$ be a smooth variety over $k$ and $Y = T^{\ast}_{X/k}$. The exterior derivative $\omega = d\theta \in \mathrm{H}^0(Y, \Lambda^2 \Omega^1_{Y/k})$ of the canonical $1$-form $\theta$ is a symplectic form.}
\begin{proof}
By definition, $\omega$ is a closed alternating $2$-form. We only need to check that it is non-degenerate. For that purpose, it is sufficient to work locally on $X$. With the identity (\ref{canonicalform}), we can write $\theta = \sum_{ i = 1}^n p_i. dx_i$ in étale coordinates, on a covering of $X$. Then, the $2$-form $\omega$ is given by: 
$$ \omega = \sum_{i =1}^n dp_i \wedge dx_i$$
and therefore non-degenerate.
\end{proof}

\subsection{Pseudo-Riemannian varieties}

{\defn A \textit{smooth pseudo-Riemannian variety over $k$} is a pair $(X,g)$ where $X$ is a smooth variety over $k$ and $g \in \mathrm{H}^0(X,S^2\Omega^1_{X/k})$ is a non-degenerate symmetric form on $X$.}

{\rem\label{symplectictangent}  Let $(X,g)$ be a smooth pseudo-Riemmanian variety over $k$. Similarly to the symplectic case, the non-degenerate symmetric form $g$ on $X$ defines an isomorphism of locally free sheaves over $X$: 
$$\Psi :\Theta_{X/k} \longrightarrow \Omega^1_{X/k} $$
such that for any local section $v$ of $\Theta_{X/k}$, $\Psi(v) = g(v,.)$. Therefore, we get an isomorphism of vector bundles over $X$:
$$T_{X/k} \simeq T^\ast_{X/k}.$$
The canonical symplectic structure on $T^\ast_{X/k}$ defined by example \ref{cotangentbundle} pulls backs via this isomorphism into a symplectic structure on $T_{X/k}$.
If $(X,g)$ is a smooth pseudo-Riemannian variety over $k$, the smooth variety $T_{X/k}$ will always be endowed with this symplectic structure.}

{\defn Let $X$ be a variety over some field $k$ and $\mathcal F$ a coherent sheaf on $X$. A \textit{connection $\nabla$ on the coherent sheaf $\mathcal F$} is a morphism of sheaves of abelian groups  $\nabla : \mathcal F \otimes \Theta_{X/k} \longrightarrow \mathcal F$, satisfying the Leibniz rule for any local sections $v,\sigma$ and $f$ of $\Theta_{X/k},\mathcal F$ and $\mathcal O_X$ respectively : 
$$\nabla_v(f.\sigma) = f. \nabla_v(\sigma) + \nabla_v(f).\sigma$$ 
where $\nabla_v(f)$ denotes the derivative of $f$ for the derivation induced by $v$.}

{\Prop\label{Levicivita} Let $(X,g)$ be a smooth pseudo-Riemannian variety over $k$. There exists a unique connection $\nabla$ on the locally free sheaf $\Theta_{X/k}$ of vector fields on $X$ such that: 
\begin{itemize}
\item[(i)] The connection $\nabla$ has zero torsion, i.e., for every vector fields $v,w \in \Theta_X(U)$ on an open set $U \subset X$, we have:
$$\nabla_v w - \nabla_w v = [v,w].$$
\item[(ii)] The metric $g$ is parallel for the connection $\nabla$, i.e., for every vector fields $v,w,x \in \Theta_X(U)$ on an open set $U \subset X$, we have:
$$\nabla_x(g(v,w)) = g(\nabla_x v, w) + g(v, \nabla_x w).$$ 
\end{itemize}}

The connection $\nabla$ on the tangent bundle of $X$ is called the \textit{Levi-Civita connection of the pseudo-Riemannian variety $X$}.

\begin{proof}
The usual formulas proving the existence and the uniqueness of the Levi-Civita connection in Riemannian geometry (see for instance \citep[Theorem 5.4]{Lee} or \citep[Part II]{Milnor}) still make sense in this context as we are in characteristic $\neq 2$. 
\end{proof}
{\defn Let $(X,g)$ be a pseudo-Riemannian variety. \textit{The curvature tensor $R$ of the pseudo-Riemannian variety $(X,g)$} is the curvature of the associated Levi-Civita connection.}

In particular, by definition, the curvature tensor $R \in \mathrm{H}^0(X, \mathrm{End}(\Theta_{X/k}) \otimes \Lambda^2 \Omega^1_{X/k})$ is an alternating $2$-form with values in the locally free sheaf $\mathrm{End}(\Theta_{X/k})$.

{\nota Let $(X,g)$ be a smooth pseudo-Riemannian variety over $k$. We denote by $\mathcal G_2(X)$ the \textit{Grassmannian of $2$-planes in $T_{X/k}$}. As $X$ is a smooth variety over $k$, the variety $\mathcal G_2(X)$ is also smooth over $k$.
We shall denote by  $\mathcal G_2^g(X)$ the open subscheme of $\mathcal G_2(X)$ whose geometric points parametrize the non-isotropic $2$-planes  (i.e., planes where the restriction of the $2$-form $g$ is non-degenerate) of $X$.}

{\cons Let $(X,g)$ be a smooth pseudo-Riemannian variety over $k$. Consider the closed subvariety  $Z \subset T_{X/k} \times_X T_{X/k}$ of pairs of orthogonal vectors in $T_{X/k}$ defined by the equation $g(v,w) = 0$ and the open subscheme $Z_0 = (U \times U) \cap Z$ of $Z$ where $U \subset T_{X/k}$ is the open subscheme defined by the inequation $g(v,v) \neq 0$.

Over any geometric point $x$ in $X$, if $(v,w) \in Z_{0,x}$, then $\mathrm{Vect}(v,w)$ is a non-isotropic $2$-plane and every isotropic $2$-plane of $T_{X/k,x}$ admits a basis in $Z_{0,x}$. Consequently, we have a surjective morphism of varieties over $k$: 
$$ \pi : Z_0 \longrightarrow \mathcal G_2^g(X).$$

Note that the curvature tensor defines a morphism $R : T_{X/k} \times_X T_{X/k} \times_X T_{X/k}  \longrightarrow T_{X/k}$ of algebraic varieties over $k$,  and that the formula 

\begin{eqnarray} \label{formulacurvature}
K^g(v,w) = \frac {g(R(v,w,v),w)}{g(v,v).g(w,w)}\end{eqnarray}
defines a regular function $K^g$ in $\mathrm{H}^0(Z_0, \mathcal O_{Z_0})$. The following lemma, together with the smoothness of $G_2^g(X)$ and the existence of sections of $\pi$ (locally in the Zariski topology) shows that this function descends to a function on $\mathcal G_2^g(X)$.

{\lem \label{sectional curvature} Let $(X,g)$ be a pseudo-Riemannian algebraic variety over $k$. Consider $\overline{k}$ an algebraic closure of $k$, $x \in X(\overline{k})$ and $P \subset T_x X(\overline{k})$ a non-isotropic $2$-plane. Consider $(v,w)$ an orthonormal basis of $P$ and set
$$K^g(P) = \frac {g(R(v,w,v),w)}{g(v,v).g(w,w)} \in \overline{k}.$$ 
The quantity $K^g(P)$ does not depend on the choice of an orthonormal basis $(v,w)$ of $P$.} 

The lemma is a simple consequence of the invariance properties   of the curvature tensor $R \in \mathrm{H}^0(X, \mathrm{End}(\Theta_{X/k}) \otimes \Lambda^2 \Omega^1_{X/k})$ (see \citep[Proposition 7.4]{Lee} for instance or \citep[Part II]{Milnor}).

{\defn\label{curvature} Let $(X,g)$ be a pseudo-Riemannian variety. We call \textit{sectional curvature} of $X$, the function $K^g \in \mathrm{H}^0(\mathcal G_2^g(X), \mathcal O_{\mathcal G_2^g(X)})$  defined by the formula (\ref{formulacurvature}) and Lemma  \ref{sectional curvature}.}

\subsection{Classical mechanics}

{\cons Let $(X,g)$ be a smooth pseudo-Riemannian variety over $k$. Consider the function $H_g \in \mathrm{H}^0(X,T_{X/k})$ associated to $g$ defined by 
$$H_g(x,p) = \frac 1 2 g_x(p,p).$$ 

By using the symplectic structure on $T_{X/k}$ given by Remark \ref{symplectictangent}, we denote by $X_g$, the Hamiltonian vector field associated to $H_g$.

{\defn\label{geodesic D-variety} Let $(X,g)$ be a smooth pseudo-Riemannian variety over $k$. We call \textit{geodesic $D$-variety associated to $(X,g)$}, the $D$-variety $(T_{X/k}, X_g)$ over the differential field $(k,0)$.}

By Lemma \ref{firstintegral}, the function $H_g$ is an integral of the $D$-variety $(T_{X/k}, X_g)$. In other words, we have a morphism of $D$-varieties over $(k,0)$:

\begin{eqnarray} \label{reduced} H_g : (T_{X/k}) \longrightarrow (\mathbb{A}^1,0). \end{eqnarray}}

Note that the only singular value of $H_g$ is $0$. Consequently, for any $t \in k \setminus \lbrace 0 \rbrace$, the fibre of $H_g$ over $t$ is a smooth variety over $k$. Moreover, the fibres of $H_g$ are invariant subvarieties of  $(T_{X/k}, X_g)$ because $H_g$ is a rational integral. This means that the vector field $X_g$ restricts to a vector field on any smooth fibre of $H_g$, which will still be denoted $X_g$. 

{\lem\label{fibre} Let $(X,g)$ be a smooth pseudo-Riemannian variety over $k$ and $a,b$ two elements of $k \setminus \lbrace 0 \rbrace$. We have an isomorphism of $D$-varieties over $(k,0)$: 
$$(H_g^{-1}(a.b^2), \frac 1 b {X_g}) \simeq (H_g^{-1}(a), X_g).$$}

\begin{proof}
This is a consequence of the homogeneity properties of the Hamiltonian $H_g$. Consider the morphism $\psi_b : T_{X/k} \longrightarrow T_{X/k}$ defined by $(x,y) \mapsto (x,b.y)$. As $b \neq 0$, the morphism $\psi_b$ is an automorphism of $T_{X/k}$.

First note that, by homogeneity, $H_g \circ \psi_b = b^2. H_g$. Therefore, $\psi_b$ restricts to an isomorphism between the subvarieties $H_g^{-1}(a)$ and $H_g^{-1}(a.b^2)$.

It remains to compute the vector field $Y = \psi_b^\ast X_g$. Let $\omega$ be the symplectic form on $T_{X/k}$. We have $\phi_b^\ast \omega = b.\omega$. The vector field $Y$ is the Hamiltonian vector field associated to the Hamiltonian $H_g \circ \psi_b$ relatively to the symplectic form $\phi_b^\ast \omega$. This means that
$$ b.\omega(- , Y) = b^2.dH = b^2.\omega(-, X_g).$$
We get that $Y = b.X_g$, which gives the isomorphism of $D$-varieties.
\end{proof}

{\rem Note that in particular, after base-changing to the algebraic closure $\overline{k}$ of $k$, all fibres (as varieties over $\overline{k}$)  of $H_g$ become isomorphic.

However, there may exist non-isomorphic fibres over $k$. For example, consider the geodesic $D$-variety of the sphere $\mathbb{S}^2 \subset \R^3$ (endowed with the Euclidean metric). In that case, the fibres over negative values and the fibres over positive values are non-isomorphic over $\R$ but become isomorphic after base-changing to $\mathbb{C}$.}

{\defn\label{reducedgeodesic D-variety} Let $(X,g)$ be a smooth pseudo-Riemannian variety over $k$. We call \textit{unitary geodesic $D$-variety associated to $(X,g)$}, the fibre (as a $D$-variety) over $\frac 1 2 \in \mathbb{A}^1(k)$ of the morphism (\ref{reduced}) of $D$-varieties over $(k,0)$.}

The underlying variety of the unitary geodesic $D$-variety is \textit{the sphere bundle $SX \subset T_{X/k}$} of $(X,g)$ defined as 
$$ SX = \lbrace (x,v) \in T_{X/k} \text{ | } g_x(v,v) = 1 \rbrace.$$

{\Prop\label{fiberchange} Let $(X,g)$ be a smooth pseudo-Riemannian variety over $k$. Suppose that the unitary  geodesic $D$-variety is absolutely irreducible and that its generic type is orthogonal to the constants.

Then, for any $c \in \overline{k} \setminus \lbrace 0 \rbrace$, the fibre of $H_g$ over $c$ (as a $D$-variety) is absolutely irreducible and its generic type is orthogonal to the constants.}  

\begin{proof}
We can suppose that $k = \overline{k}$ as everything is invariant under base-change.
Fix $c \in \overline{k} \setminus \lbrace 0 \rbrace$ and consider $b$ a square root of $2.c$. By Lemma \ref{fibre}, we get an isomorphism of $D$-varieties over $k$

$$(H_g^{-1}(c), \frac 1 b {X_g}) \simeq (H_g^{-1}(\frac 1 2), X_g).$$

By hypothesis, the latter $D$-variety is irreducible and its generic type is orthogonal to the constants. Therefore, the $D$-variety $(H_g^{-1}(c), \frac 1 b {X_g})$ is also irreducible and its generic type is orthogonal to the constants.

Note that for all $n \in \mathbb{N}$, the rational first integrals of $(H_g^{-1}(c), \frac 1 b {X_g})^n$ and $(H_g^{-1}(c), {X_g})^n$ are the same. Using Theorem 2.3.4 in \cite{moi}, we conclude that the generic type of $(H_g^{-1}(c), {X_g})$ (which is the fibre of $H_g$ over $c$) is orthogonal to the constants.  
\end{proof}

{\rem In classical mechanics, the Hamiltonian $H_g(x,p) = \frac 1 2 g_x(p,p)$ can be interpreted as the kinetic energy of a particle of mass $m = 1$ moving on the pseudo-Riemannian manifold $(X,g)$.
Therefore, the geodesic $D$-variety describes the movement of a free particle moving on the pseudo-Riemannian manifold $(X,g)$.}   

More generally, the preceding ideas are also relevant to the case of a particle subject to a potential $V \in \mathrm{H}^0(X, \mathcal O_X)$.

{\defn Let $(X,g)$ be a smooth pseudo-Riemannian variety over $k$ and $V \in \mathrm{H}^0(X,\mathcal O_X)$ a function on $X$. The $D$-variety \textit{$(T_{X/k}, X_{V,g})$  of the pseudo-Riemannian variety $(X,g)$ with potential $V$} is the $D$-variety over $(k,0)$ defined by the Hamiltonian $H_{g,V} \in \mathrm{H}^0(X,T_{X/k})$ :
$$H_{g,V}(x,p) = \frac 1 2 g_x(p,p) + V(x).$$}

{\exam[$n$-body problem] Set $M = \R^3$ and consider the Riemannian manifold $(M^n \setminus \Delta, g_{0})$ where $\Delta$ is the union of the diagonals and $g_0$ is the restriction of the Euclidean metric to $M^n \setminus \Delta$. 

Recall that, in classical mechanics, the $n$-body problem is the  Hamiltonian system for the Riemannian manifold $(M^n \setminus \Delta, g_{0})$ associated to the potential 
$$V(q_1,\cdots, q_n) = - \sum_{1 \leq i < j \leq n} \mathcal G \frac {m_i.m_j}{|| q_i - q_j ||}.$$}

where $\mathcal G$ is the gravitational constant and the $m_i$ are the masses of the bodies (see \cite{Chen}). Notice that the potential $V$ is an analytic (non polynomial) function on $M^n \setminus \Delta$.

{\cons Consider the subset of $(M^n  \setminus \Delta) \times \R^{\frac {n(n-1)} 2}$ defined by 
$$ \tilde{M} := \lbrace ((q_i)_{1 \leq i \leq n}, (z_{k,l})_{1 \leq k < l \leq n}) \in (M^n \setminus \Delta)\times \R^{\frac {n(n-1)} 2} \text{ | } z_{k,l}^2 = ||q_i - q_j||^2 \rbrace.$$

It is easy to see that the projection on $(M^n  \setminus \Delta)$ identifies it with an étale covering of $(M^n  \setminus \Delta)$ of degree $2^{\frac {n(n-1)} 2}$. Moreover, note that the function
$$||-||^2 : \begin{cases} M^2 \setminus \Delta \longrightarrow \R \\
(q_1,q_2) \mapsto  ||q_1 - q_2||^2 = (q_1^1 - q_2^1)^2 + (q_1^2 - q_2^2)^2 + (q_1^3 - q_2^3)^2 \end{cases}$$
is polynomial. Therefore, there exists a closed subvariety  $\tilde{X}$ of $((\mathbb{A}^3)^n \setminus \Delta) \times \mathbb{A}^{\frac {n(n-1)} 2}$ such that $\tilde{X}(\R) = \tilde{M}$.}

{\defn We call \textit{algebraic $n$-body problem}, the $D$-variety of the pseudo-Riemannian variety $(\tilde{X},\tilde{g_0})$ associated with the potential
$$V = -  \sum_{1 \leq i < j \leq n} \mathcal G \frac {m_i.m_j}{z_{i,j}}$$
where $\tilde{g_0}$ is the pull-back to $\tilde{X}$ of the metric $g_0$ on $\mathbb{A}^n  \setminus \Delta$.}

\section{Approximating compact Riemannian manifolds by pseudo-Riemannian varieties}

In this section, we prove the approximation results, Theorem \ref{Int1} and \ref{Int2}, stated in the introduction.
We start by recalling the notions from real algebraic geometry needed for the proof of Theorem \ref{Int1} and \ref{Int2}, before analysing the properties of the real analytification functor. The definition and results that we use are based on the chapters 12 and 14 in \cite{Coste}.

Note that in this section, there are three different categories of geometric objects (namely algebraic varieties, real-analytic varieties and smooth manifolds) and that this requires some care with the terminology. In particular, the adjective smooth will be used for objects (functions, morphisms and manifolds) from classical real differential geometry. 
\subsection{Real algebraic sets and analytification}

{\defn Let $X$ be an algebraic variety over $\R$. The \textit{real algebraic set associated to $X$} is the locally ringed space $(M = X(\R), \mathcal R_M)$ where
\begin{itemize}
\item $X(\R)$ is the set of real points of $X$ endowed with the Zariski topology;
\item the sheaf of rings $\mathcal R_M$ is the sheaf associated with the presheaf $\mathcal O_{X,\R}$ induced by $\mathcal O_X$ on $X(\R)$ and called the \textit{sheaf of regular functions}.
\end{itemize}}

Let $(X,\mathcal O_X)$ be an irreducible variety over $\R$ and $(M,\mathcal R_M)$ be the associated real algebraic set. By construction, we have a morphism of $\R$-algebras: 
$$\mathcal O_X(X){\longrightarrow} \mathcal R_M(M) $$
which is injective if $X(\R)$ is Zariski-dense in $X$. \textit{The ring of polynomial functions on $M$} is the image of this morphism. 

{\rem It is possible to describe the ring $\mathcal R_M(M)$ of regular functions on $M$ as the subring of $\mathbb{R}(X)$ of rational functions everywhere defined on $X(\R)$.
In particular as soon as $\mathrm{dim}(X) \geq1$, the $\R$-algebra $\mathcal R_M(M)$ is not of finite type.

In the sequel, we will be mainly interested in polynomial functions (and more generally polynomial sections of coherent sheaves) on $M$ although one has to introduce $\mathcal R_M$ in order to have a sheaf. Note that if $X \subset \mathbb{A}^n$ is an affine closed subvariety, then the polynomial functions on $M$ are simply the restrictions to $M$ of the polynomial functions on $\R^n$.}

{\defn Let $X$ be a variety over $\R$ and $x \in M = X(\R)$. We say that \textit{the real algebraic set associated to $X$ is non-singular at $x \in X$} if $\mathcal O_{X,x} = \mathcal R_{M,x}$ is a regular local ring.}

The real algebraic set associated to $X$ is said to be \textit{non-singular} if it is non singular at any of its points, or equivalently if $X(\R)$ is contained in the open set $X_{reg}$ of regular points of $X$. 

{\rem There exists a \textit{real analytification functor} from the category of algebraic varieties over $\R$ to the category of real analytic spaces. If $X$ is a variety over $\R$ and $M = X(\R)$ is non-singular then $M^{an}$ is an analytic (so in particular smooth) manifold such that:
\begin{itemize}
\item The topological space $M^{an}$ is the set $X(\R)$ endowed with the Euclidean topology.
\item Every regular function on $M$ is smooth on $M^{an}$.
\end{itemize}}

{\defn Let $X$ be an algebraic variety over $\R$. Suppose that the associated real algebraic set $M = X(\R)$ is non-singular. For any coherent sheaf $\mathcal F$ on $X$, the \textit{analytification of $\mathcal F$} denoted $\mathcal F^{an}$ is the $\mathcal C^\infty(M^{an})$-module defined by: for any $x \in M$,
$$ \mathcal F^{an}_x = \mathcal F_{x,\R} \otimes _{\mathcal R_{M,x}} \mathcal C^\infty(M^{an})_x.$$}

In particular, if $\mathcal F$ is a coherent sheaf, then we have a morphism of $\R$-vector spaces
$$ -_\R^{an} : \mathrm H^0(X,\mathcal F) \longrightarrow \mathrm H^0(M,\mathcal F_\R) \longrightarrow \mathcal C^\infty(M^{an},\mathcal F^{an})$$ which is injective if $X(\R)$ is Zariski-dense in $X$.

{\rem If $\mathcal F$ is a locally free sheaf on $X$ then $\mathcal F^{an}$ is also a locally free sheaf on $M = X(\R)^{an}$. Moreover, if $\mathcal F$ is the sheaf of tensors of type $(n,m)$, then so is $\mathcal F^{an}$.}

Recall that \textit{a pseudo-Riemannian (analytic) manifold} $(M,g)$ is an analytic manifold $M$ endowed with a non-degenerate symmetric 2-form $g$. The pseudo-Riemannian manifold $(M,g)$ is called \textit{a Riemannian manifold} if the $2$-form $g$ is positive definite or, in other words, if it has signature $(n,0)$.

{\defn Let $(X,g)$ be a smooth pseudo-Riemannian variety over $\mathbb{R}$. The \textit{analytification of $(X,g)$} is the pseudo-Riemannian (analytic) manifold $(X(\R)^{an}, g_\R^{an})$. It is denoted $(X(\R),g_\R)^{an}.$}

The following two lemmas relate the notions of the first section of this article on pseudo-Riemannian varieties with the classical ones on their analytification.

{\lem Let $(X,g)$ be a smooth pseudo-Riemannian variety over $\R$. Denote by $K^g$ the sectional curvature of $(X,g)$. 

The sectional  curvature (in the sense of differential geometry) of the pseudo-Riemannian manifold $(X(\R),g_\R)^{an}$ is the analytification $(K_g)^{an}_\R : \mathcal G^g_2(X(\R)^{an}) \longrightarrow \R$ of $K_g$.}

\begin{proof}
By uniqueness of the Levi-Civita connection in Riemannian geometry, the Levi-Civita of $(X(\R),g_\R)^{an}$ is $\nabla_\R^{an}$ where $\nabla$ is the Levi-Civita connection of the smooth pseudo-Riemannian variety $(X,g)$. Therefore, if $R$ denotes the curvature tensor of $\nabla$, the curvature tensor of $\nabla_\R^{an}$ is $R_\R^{an}$. Then, the sectional curvature is computed from $R$ (Lemma \ref{sectional curvature}) with the same formula as in Riemannian geometry.
\end{proof}

{\lem\label{realflow} Let $(X,g)$ be a smooth pseudo-Riemannian manifold over $\R$.  The real flow of the unitary $D$-variety (resp. the reduced geodesic $D$-variety) associated to $(X,g)$ is the geodesic flow (resp. unitary geodesic flow) of the Riemannian manifold $(X(\R),g_\R)^{an}$.}

\begin{proof}
Set $M = X(\R)^{an}$. It is easy to see that $T^\ast_M = T^\ast_{X/\R}(\R)^{an}$ and that the canonical symplectic structure on $T^\ast_M$ is given by the closed non-degenerate $2$-form $\omega_\R^{an}$ where $\omega$ is the symplectic form $T^\ast_{X/\R}$ given by example \ref{cotangentbundle}.

Moreover, the isomorphism $T_{X/\R} \simeq T^{\ast}_{X/\R}$ of varieties over $k$ given by $g$ induces the isomorphism $TM \simeq T^\ast M$ defined by $g_\R^{an}$. Consequently, if $H$ is a Hamiltonian on $T_{X/\R}$, the Hamiltonian vector field on $TM$ associated to $H_\R^{an}$ is $(X_{H})_\R^{an}$.
 
By definition, the geodesic flow is the Hamiltonian flow on $TM$ associated with the Hamiltonian $H_{g_\R^{an}}(x,p) = g_x^{an}(p,p)$. As $H_{g_\R^{an}} = H_g^{an}$, the lemma follows from the previous compatibilities. 
\end{proof}

\subsection{First approximation} We prove that every compact Riemannian manifold can be approximated  by a smooth pseudo-Riemannian manifold over $\R$. 

{\nota Let $M$ be a smooth and compact manifold of dimension $n \in \mathbb{N}$. Recall that the ring $\mathcal C^\infty(M)$ of smooth functions on $M$ has a natural topology of a Fréchet space, called the \textit{$\mathcal C^\infty$--topology} defined by the family of seminorms
$$||f||_{N,U} = \mathrm{sup}_{x \in U} \lbrace \partial^\alpha f(x) \text{ | } |\alpha | \leq N \rbrace.$$
where $U \subset_i M$ is a chart of the manifold $M$.

More generally, if $E \longrightarrow M$ is a vector bundle of rank $r$, the vector space $\mathcal C^\infty(E)$ of smooth sections of $E$ also carries the natural topology of a Fréchet space called the \textit{$\mathcal C^\infty$--topology} defined as follows.

Consider any trivialising covering $(U_i)_{i = 1}^n$ of $M$. This covering defines an inclusion 
$$ \mathcal C^\infty(E) \longrightarrow \Pi_{i = 1}^n \mathcal [C^\infty(U_i)]^r.$$

The smooth topology on $\mathcal C^\infty(E)$ is the topology induced by the product topology on $\Pi_{i = 1}^n \mathcal [C^\infty(U_i)]^r$. It does not depend on the choice of the trivialising covering.}
 
{\Thm[Stone-Weierstrass Theorem -- {\citep[Théorème 8.8.5]{Coste}}]\label{Stone-Weierstrass} Let $M \subset \R^n$ be a smooth and compact submanifold. The restrictions to $M$ of polynomial functions on $\R^n$ are dense in $\mathcal C^\infty(M)$ endowed with the $\mathcal C^\infty$--topology.}

{\cor\label{smoothsheaf} Let $X$ be an affine variety over $\R$ and $\mathcal F$ a locally free sheaf on $X$. Suppose that the real algebraic set $M = X(\R)$ associated to $X$ is non-singular and compact. The morphism of  $\R$-vector spaces:
$$-_\R^{an} : \mathrm{H}^0(X, \mathcal F) \longrightarrow \mathcal C^\infty(M^{an},\F^{an}) $$
has a dense image for the $\mathcal C^\infty$--topology.}

\begin{proof}

We fix a closed embedding $X \subset_i \mathbb{A}^n$. The analytic morphism $i_{\R}^{an}$ realises $M^{an}$ as a compact smooth submanifold of the Euclidean space $\R^n$. Therefore, Theorem \ref{Stone-Weierstrass} above shows that the corollary is true when $\mathcal F = \mathcal O_X$ and therefore when  $\mathcal F = \mathcal O_X^n$ is a trivial bundle.

Let $\mathcal F$ be a locally free sheaf on $X$. As $X$ is an affine variety, the coherent sheaf $\mathcal F$ is generated by its global sections. Therefore, there exists a surjective morphism: 
$$\phi : \mathcal O_X^n \longrightarrow \mathcal F.$$

The morphism $\phi_\R^{an} : \mathcal O_M^n \longrightarrow \mathcal F^{an}$ of analytic sheaves is also surjective. Consequently, to approximate a smooth section $\sigma$ of $\mathcal F^{an}$, it is sufficient to lift it up into a section $\tau$ of $\mathcal O_M^n$ and then to approximate $\tau$.
\end{proof}

{\Thm[Tognoli's Theorem --{\citep[Théorème 14.1.10]{Coste}}]\label{Tognoli} Let $M$ be a smooth, connected and compact manifold. There exists an affine variety $X$ over $\R$ such that $X(\R)$ is non-singular and the smooth manifold $X(\R)^{an}$ is isomorphic to $M$.}

{\rem\label{Nashtheorem} A weak version of Tognoli's Theorem was already established by Nash. It states that any connected and compact manifold is a connected component of an algebraic subset of $\R^n$ for some $n \in \mathbb{N}$.

In fact, the proof of Nash's result is quite elementary.  Fix a smooth embedding $i : M \longrightarrow \R^n$ and consider a tubular neighbourhood $T$ of $M$ in $\R^n$. By construction, there exist smooth functions $f_1,\ldots, f_r$ on $T$ such that 
$$ M = \lbrace x \in T \text{ | } f_i(x) = 0 \text{ , } \forall i \leq n \rbrace.$$

By the Thom Transversality Lemma (see \citep[Section 14]{Coste}), for any sufficiently small pertubations $(g_i)_{1 \leq i \leq n}$ of the $(f_i)_{1 \leq i \leq n}$ on $\overline{T}$, the resulting smooth manifold  $ \tilde{M} = \lbrace x \in T \text{ | } g_i(x) = 0 \text{ , } \forall i \leq n \rbrace$ will be diffeomorphic to $M$ and still satisfy $\tilde{M} \cap \partial \overline{T} = \emptyset$.

By Theorem \ref{Stone-Weierstrass}, we can choose the $(g_i)_{i \leq n}$ to be the restrictions to $\overline{T}$ of polynomials $(P_i)_{i \leq n}$. As $\tilde{M} \cap \partial \overline{T} = \emptyset$, we can conclude that $\tilde{M} = \overline{T} \cap Z(P_1,\cdots, P_n)$ is both closed and open in $Z(P_1,\cdots, P_n)$ and is therefore a connected component.

Consequently, the hard part in Tognoli's Theorem is to get rid of the other connected components.}

{\nota If $M$ is a smooth manifold then $\mathcal C^\infty(S^2T^\ast M)$ denotes the vector space of global $2$-forms on $M$ which are the sections of the vector bundle  $S^2(T^\ast M)$.
The Riemannian metrics on $M$ form an open subset of $\mathcal C^\infty(S^2T^\ast M)$.}

{\Thm\label{firstapproximation} Let $(M,g_M)$ be a smooth, connected and compact Riemannian manifold. For every neighbourhood $\mathcal U \subset \mathcal C^\infty(S^2T^\ast M)$ of $g_M$ for the $\mathcal C^\infty$--topology, there exist a smooth pseudo-Riemannian variety $(X,g)$ over $\R$ and a diffeomorphism $\phi : M \longrightarrow X(\R)^{an}$ such that $\phi^\ast g^{an}_{\R} \in \mathcal U$.} 

\begin{proof}
Let $(M,g_M)$ be a smooth, connected and compact Riemannian manifold. By Theorem  \ref{Tognoli}, there exists an affine variety $X$ over $\R$ such that the real algebraic set $X(\R)$ is non-singular and $M$ and $X(\R)^{an}$ are diffeomorphic.

Fix a neighbourhood $\mathcal U \subset \mathcal C^\infty(S^2M)$ of $g_M$. Up to refining $\mathcal U$, we can suppose that any $g \in \mathcal U$ is non-degenerate. By corollary \ref{smoothsheaf} applied to the affine variety $X$, there exists a section $g \in \mathrm{H}^0(X,S^2(T_{X/\R}))$ such that $\phi^\ast g^{an}_{\R} \in \mathcal U$.

Let $V \subset X$ be the biggest open set such that: 
\begin{itemize}
\item All points of $V$ are regular points.
\item The symmetric $2$-form $g$ is non degenerate.
\end{itemize}
By construction, $(V,g_{|V})$ is a smooth pseudo-Riemannian variety over $\R$. 
Moreover, as $X(\R)$ is non-singular and $g_\R^{an}$ (and therefore $g$) is non-degenerate on it, we have $X(\R) = V(\R)$. Consequently, $(V,g_{|V})$ and $\phi : M \longrightarrow V(\R) = X(\R)$ satisfy the conclusion of the theorem. 
\end{proof}

\subsection{Approximation by algebraic embeddings in $\R^n$} 

{\Thm\label{secondapproximation} Let $(M,g)$ be a smooth, connected and compact Riemannian manifold. For every neighbourhood $\mathcal U \subset \mathcal C^\infty(S^2T^\ast M)$ of $g_M$ for the $\mathcal C^\infty$--topology, there exist $m \in \mathbb{N}$ and a $\mathcal C^\infty$-embedding $i :  M \longrightarrow \R^m$ such that :
\begin{itemize}
\item[(i)] The image of $M$ is a non-singular algebraic subset $X(\R)$ of $\R^m$, defined by some smooth, locally closed subvariety $X$ of $\mathbb{A}^m_\R$.
\item[(ii)] Let $g_{0}$ be the usual Euclidean metric on $\R^m$. The pullback $i^\ast g_0$ of the Riemannian metric $g_0$ is in $\mathcal U$. 
\end{itemize}}

The proof is based on a clever elementary lemma that goes back to Nash in \cite{Nash2} and \cite{Nash} (see also \citep[Section 3.1]{Gromov} for a modern presentation).

{\lem[Nash Twist]\label{nashtwist} Let $M$ be a compact manifold and $f,\phi : M \longrightarrow \R$ be smooth functions. Fix $\epsilon > 0$, then there exist smooth functions $h : M \longrightarrow \R$  and $k : M \longrightarrow \R$ such that:  
$$\phi^2 (df)^2 = (dh)^2 + (dk)^2 - \epsilon^2 (d\phi)^2.$$}

\begin{proof}
Consider the function $f_0$ from $M$ to the circle $\mathbb{S}_\epsilon^1 \subset \R^2$ of radius $\epsilon$ defined by:
$$f_0(p) = (h_0(p),k_0(p)) = (\epsilon.\mathrm{cos}(\epsilon^{-1} f(p)) , \epsilon. \mathrm{sin}(\epsilon^{-1} g(p)).$$
Note that $f_0$ as been chosen such that the pull-back of the Euclidean metric on $\R^2$ by $f_0$ is $(df)^2$.

Set $f_1 = \phi.f_0$ and $g_1$ the pullback by $f_1$ of the Euclidean metric on $\R^2$. Using the Leibniz rule, we can compute:
$$g_1 = \phi^2.(df)^2 + \epsilon^2. (d\phi)^2.$$
Consequently, if $h$ and $k$ are the coordinate functions of $f_1$, we have $(dk)^2 + (dh)^2 = \phi^2.(df)^2 + \epsilon^2. (d\phi)^2$.
\end{proof}

{\Prop \label{keyapproximationlemma} Let $X$ be an affine variety over $\R$ such that $X(\R) = M \subset \R^N$ is a compact and non-singular algebraic subset. Consider  smooth functions $\phi,f : X(\R)^{an} \longrightarrow \R$  on $X(\R)^{an}$. 

For any neighbourhood $\mathcal U \subset \mathcal C^\infty(S^2 T^\ast M)$ of $0$, there  exist $h,k \in \R[X]$ such that: 
$$\phi^2 df^2 - (dh^2 + dk^2) \in \mathcal U.$$}

\begin{proof}
Let $X$ be an affine variety over $\R$ such that $X(\R) = M \subset \R^N$ is a compact and non-singular algebraic subset.

Consider smooth functions $\phi,f$ on $X(\R)^{an}$ and fix $\mathcal U \subset \mathcal C^\infty(S^2 T^\ast M)$ a neighbourhood of $\phi^2 (df)^2$.
Let $\epsilon > 0$ be such that $- \epsilon^2 (d\phi)^2 \in \mathcal U$ and let $h$ and $k$ be the smooth functions on $M$ given by Lemma \ref{nashtwist} such that: 
\begin{eqnarray}\label{nahformula}
\phi^2 df = (dh)^2 + (dk)^2 - \epsilon^2 (d\phi)^2.
\end{eqnarray}

By Theorem \ref{Stone-Weierstrass}, there exist sequences of polynomial functions $h_n \in \R[X]$ and $k_n \in \R[X]$ converging to $h \in \mathcal C^\infty(M)$ and $k \in \mathcal C^\infty(M)$ respectively for the $\mathcal C^\infty$--topology.  

Using the equation (\ref{nahformula}), we have: 
$$ \phi^2.df - [(dh_n)^2 + (dk_n)^2] = [(dh)^2 - (dh_n)^2]  + [(dk)^2  - (dk_n)^2] - \epsilon^2 (d\phi)^2 \overset{n \mapsto \infty}{\longrightarrow} - \epsilon^2 (d\phi)^2$$
As $- \epsilon^2 (d\phi)^2 \in \mathcal U$ and $\mathcal U$ is open, we conclude that for $N$ large enough we have: 
$$ \phi^2.df - [(dh_N)^2 + (dk_N)^2] \in \mathcal U. \qedhere$$
\end{proof}

\begin{proof}[Proof of Theorem \ref{secondapproximation}]
Let $(M,g)$ be a smooth, connected and compact Riemannian manifold of dimension $n \in \mathbb{N}$. We start by a couple of reductions.

By Theorem \ref{Tognoli}, there exists an affine variety $X$ over $\R$ such that $X(\R)$ is non-singular and $X(\R)^{an} \simeq M$. Consequently, we can suppose that $M = X(\R)^{an} \subset \R^N$ is the analytification of an algebraic subset of $\R^N$.

Let $g_0$ be the restriction of the Euclidean metric of $\R^N$ to the submanifold $M$. Then $g - g_0$ is a symmetric $2$-form which might not be positive definite.
But, as $X(\R)$ is compact, by replacing $X$ by its image by a homothety (which is an algebraic automorphism), we can suppose that $||g_0||_\infty$ is sufficiently small for $g - g_0$ to be a positive definite $2$-form on $M$. 
Consequently, there exist smooth functions $(\phi_i)_{i \leq L}$ and $(\psi_i)_{i \leq L}$  on $M$ such that:
$$g - g_0 = \sum_{i = 1}^L \phi_i^2.(d\psi_i)^2 $$
(see for instance \cite[pp. 222]{Gromov} for the existence of such a decomposition).
Fix a neighbourhood $\mathcal U \subset \mathcal C^\infty(S^2 T^\ast M)$ of $g$ and consider for each $i \leq L$, a neighbourhood $\mathcal U_i \subset \mathcal C^\infty(S^2 T^\ast M)$ of $0$ such that
\begin{eqnarray}\label{decomposition metric}
[\mathcal U_1 + \phi_1^2(d\psi_1)^2] + \cdots + [\mathcal U_N + \phi_N^2(d\psi_L)^2] \subset \mathcal U + g_0.
\end{eqnarray}

By Proposition \ref{keyapproximationlemma}, for all $i \leq L$, there exists a morphism $\rho_i = (h_i,k_i) : X \longrightarrow \mathbb{A}^2$ satisfying 
\begin{eqnarray}\label{Gromov} \phi_i^2 (d\psi_i)^2 - (dh_i^2 + dk_i^2) \in \mathcal U_i \end{eqnarray}

Now, consider the morphism $\rho = (\rho_1,\cdots , \rho_L) : X \longrightarrow (\mathbb{A}^2)^L$ and $Y = \mathrm{graph}(\rho) \subset \mathbb{A}^N \times (\mathbb{A}^2)^L$. First note that the morphism  $\chi : X \longrightarrow Y$ defined by $x \mapsto (x,\rho(x))$ is an isomorphism of varieties over $\R$. Therefore, the metric $(\chi_\R^{an})^\ast g_1$ is given as the restriction of the Euclidean metric to an algebraic embedding of $M$. 

Moreover, by construction, the restriction $g_1$ of the Euclidean metric to $Y$ is given by: 
$$ (\chi_\R^{an})^\ast g_1 = g_0 + \sum_{i = 1}^L (dh_i^2 + dk_i^2).$$

By using the decomposition (\ref{decomposition metric}) and the property (\ref{Gromov}), we get that $(\chi_\R^{an})^\ast g_1 - g \in \mathcal U$.
\end{proof}

\section{Anosov flows and their dynamics}
In this section, we give a self-contained exposition of the mixing properties of the geodesic flow of a compact Riemannian manifold with negative curvature, involved in the proof of Theorem \ref{Int3}.

First, we recall how to translate the hypothesis on the curvature of the Riemannian manifold into a global hyperbolic structure for the geodesic flow -- \textit{the structure of an Anosov flow} -- discovered by D. Anosov in \cite{Ano}.

Then, we gather the consequences of the existence of this hyperbolic structure on the topological dynamics of the geodesic flow. Those results will allow us to deduce some mixing properties for the geodesic flow from statements concerning the length of the periodic orbits and the density of recurrent points of the flow (see for instance, \cite{Coud}).

It is then easy to deduce the mixing properties of the geodesic flow for compact Riemannian manifold with negative curvature, in the specific form needed for the proof of Theorem \ref{Int3}, from the work of F. Dal'bo in \cite{Dalbo}.

\subsection{Anosov Flows}

{\defn Let $(M,g)$ be a Riemannian manifold of dimension $\geq 3$ and $v$ a $\mathcal C^\infty$-vector field on $M$ with a complete flow $(\phi_t)_{t \in \R}$. The flow $(M , (\phi_t)_{t \in \R})$ is called \textit{an Anosov flow} if the vector field $v$ does not vanish and there exists a splitting of the tangent bundle $TM$ into continuous sub-bundles
\begin{eqnarray}\label{Anosovsplitting}
TM = E^s \oplus \R.v \oplus E^u \end{eqnarray}
satisfying:
\begin{itemize}
\item[(i)] The sub-bundles $E^s$ and $E^u$ are non-trivial bundles which are $(d\phi_t)_{t \in \R}$-invariant.
\item[(ii)] There exist $C,C' > 0$ and $0 < \lambda < 1$ such that for all $u \in E^s$,
$$|| d\phi_t(u) || \leq  C. \lambda^t ||u|| \text{ and } || d\phi_{-t}(u) || \geq  C'. \lambda^{-t} ||u|| \text{ for all t > 0}.$$
\item[(ii)] There exist $C,C' > 0$ and $0 < \lambda < 1$ such that for all $w \in E^u$,
$$|| d\phi_t(w) || \geq  C. \lambda^{-t} ||w|| \text{ and } || d\phi_{-t}(w) || \leq  C'. \lambda^{t} ||w|| \text{ for all t > 0}.$$
where the norm on the tangent bundle $TM$ is given by the Riemannian metric on $M$. 
\end{itemize}}

Let $M$ be a compact manifold. If $g_1$ and $g_2$ are two Riemannian metrics on $M$, then the associated norms $||.||_1$ and $||.||_2$ are equivalent, in the sense that there exist $A,B > 0$ such that 
$$A. ||.||_1 \leq ||.||_2 \leq  B.||.||_1$$

In particular, if $v$ is a vector field on a compact manifold $M$, then the property of being Anosov does not depend on the choice of the Riemannian metric on $M$. Accordingly, we will speak about Anosov flow on a compact manifold $M$ without a specific choice of Riemannian metric on $M$.

\rem Let $M$ be a Riemannian manifold and $v$ a vector field on $M$ with a complete flow $(\phi_t)_{t \in \R}$. For every $x \in M$, one can always define the two vector subspaces of $T_x M$ by: 
$$ E_x^s = \lbrace v \in T_x M \text{ | } || d\phi_t(v) || \overset{t \mapsto \infty}{\longrightarrow} 0 \rbrace \text{ and }
 E_u^s = \lbrace v \in T_x M \text{ | } || d\phi_t(v) || \overset{t \mapsto - \infty}{\longrightarrow} 0 \rbrace.$$
The vector subspaces $E^s_x$ and $E^u_x$ of $T_x M$ are respectively called the \textit{(strongly) stable} and the \textit{(strongly) unstable distributions at $x \in M$} of the flow $(\phi_t)_{t \in \R}$.

Note that, for a general vector field $v$, the collection $(E^s_x)_{x \in M}$ does not behave continuously, when $x$ varies in $M$. However, when the flow of $v$ is an Anosov flow, it is easy to check that the stable and unstable distributions at $x \in M$ coïncide with the continuous subbundles defining the splitting (\ref{Anosovsplitting}) in the definition of an Anosov flow. Consequently, when the flow of $v$ is an Anosov flow, the collection $(E^s_x)_{x \in M}$ (and similarly for the unstable distribution) defines a continuous vector subbundle of $TM$.

Moreover, a consequence of the previous analysis is that the splitting $(\ref{Anosovsplitting})$ is unique, if it exists. We can therefore rephrase the definition of an Anosov flow in the following terms. The complete flow $(\phi_t)_{t \in \R}$ associated to a vector field $v$ is an Anosov flow  if:
\begin{itemize}
\item The stable and unstable distributions are both (non-trivial) continuous vector bundles and define a splitting $TM = E^s \oplus \R.v \oplus E^u$ of the tangent bundle of $M$.
\item We have exponential bounds for the rate of convergence and divergence in the stable and unstable distributions $E^s$ and $E^u$.
\end{itemize}

Regarding the unitary geodesic flow on a compact Riemannian manifold with negative curvature, we have the following theorem. 

{\Thm[{\cite{Eber}}]\label{anosov} Let $(M,g)$ be a compact Riemannian manifold with negative curvature. The unitary geodesic flow $(SM,(g_t)_{t \in \R})$ on the sphere bundle of $M$ is an Anosov flow.} 

\subsection{Dynamical properties of Anosov flows} We are now interested in the \textit{ topological dynamics} of an Anosov flow. The natural setting for this section is therefore the one of (complete) \textit{ metric flows}, i.e., of metric spaces $(M,d)$ endowed with a complete flow $(\phi_t)_{t \in \R}$.

{\defn Let $(M,(\phi_t)_{t \in \R})$ be a metric flow and $p \in M$. The \textit{(strongly) stable set} at $p$, denoted by $W^s(p)$, and \textit{the (strongly) unstable set} at $p$, denoted by $W^u(p)$, are the subsets of $M$ defined by:  
\begin{eqnarray*}
W^s(p)  = & \lbrace q \in M \text{ | } d(\phi_t(p),\phi_t(q)) \overset{t \mapsto \infty} \longrightarrow 0 \rbrace \\
W^u(p) = & \lbrace q \in M \text{ | } d(\phi_{-t}(p),\phi_{-t}(q)) \overset{t \mapsto  \infty} \longrightarrow 0 \rbrace.
\end{eqnarray*}}

{\defn Let $(M,(\phi_t)_{t \in \R})$ be a metric flow and $p \in M$. For any $\epsilon > 0$, we call \textit{strongly stable set at $p$ of size $\epsilon$} and \textit{strongly unstable set at $p$ of size $\epsilon$}, the sets defined by:
\begin{eqnarray*}
W_\epsilon^s(p)  =& \lbrace q \in W^s(p) \text{ | } d(\phi_t(p),\phi_t(q)) \leq \epsilon \text{ for all } t \geq 0 \rbrace \\
W_\epsilon^u(p)  =& \lbrace q \in W^u(p) \text{ | } d(\phi_{-t}(p),\phi_{-t}(q)) \leq \epsilon \text{ for all } t \geq 0 \rbrace.
\end{eqnarray*}}
  
Let $(M,g)$ be a Riemannian manifold and $v$ a $\mathcal{C}^\infty$-vector field on $M$ with a complete flow $(\phi_t)_{t \in \R}$. When $(M,(\phi_t)_{t \in \R})$ is an Anosov flow then, for every $\epsilon > 0$, the germ at $p$ of the strongly stable set at $p$ of size $\epsilon$ is the germ at $p$ of the strongly stable set. Namely we have:

{\lem[{\citep[Section 2.2.i]{Hass}}] Let $(M,(\phi_t)_{t \in \R})$ be an Anosov flow as above. For every $\epsilon > 0$, there exists $\delta >0$ such that:
$$W^s(p) \cap B(p,\delta) = W_\epsilon^s(p) \cap B(p,\delta).$$}

{\Thm[Hadamard-Perron Theorem -- {\citep[Section 2.2.i]{Hass}}] Let $(M,(\phi_t)_{t \in \R})$ be an Anosov flow as above. The stable and unstable distributions $E^s$ and $E^u$ are integrable by foliations with $\mathcal C^1$ leaves. Moreover, for all $p \in M$, the leaf at $p \in M$ of the associated foliation to $E^s \subset TM$ (resp. to $E^u \subset TM$) is $W^s(p)$ (resp. $W^u(p)$).}

{\defn Let $(M,(\phi_t)_{t \in \R})$ be a metric flow. We say that the flow 
$(M,(\phi_t)_{t \in \R})$ admits a \textit{local product structure} if for all $p \in M$, there exists a neighbourhood $V \subset M$ of $p$ which satisfies:

For every $\epsilon > 0$, there is a positive constant $\delta > 0$ such that for all $x,y \in V$ satisfying $d(x,y) \leq \delta$, there exists $t \in \R$ with $|t| \leq \epsilon$ such that 
$$ W^u_\epsilon (\phi_t(x)) \cap W^s_\epsilon (y) \neq \emptyset.$$}

{\Prop[{\cite{Ano}}]\label{localproductstructure} Let $(M,(\phi_t)_{t \in \R})$ be an Anosov flow as above. Then $(M,(\phi_t)_{t \in \R})$ admits a local product structure.} 

{\defn Let $(M,(\phi_t)_{t \in \R})$ be a metric flow. We say that the flow 
$(M,(\phi_t)_{t \in \R})$ satisfies the  \textit{Anosov closing lemma} if for all $p \in M$, there exists a neighbourhood $V \subset M$ of $p$ which satisfies:

For all $\epsilon > 0$, there exist positive constants $\delta > 0$ and $L > 0$ such that for all $x \in V$ and $t > L$ with $d(x,\phi_t(x)) \leq \delta$ and $\phi_t(x) \in V$, there exist $l > 0$ and  $x_0 \in V$ satisfying: 
\begin{itemize}
\item[(i)] $|l - t| \leq \epsilon$ and $\phi_l(x_0) = x_0$;
\item[(ii)] For all $s \leq \mathrm{min}(t,l)$, $d(\phi_s(x_0), \phi_s(x)) \leq \epsilon$.
\end{itemize}

{\Prop[{\cite{Ano}}]\label{Anosovclosinglemma} Let $(M,(\phi_t)_{t \in \R})$ be an Anosov flow as above. Then $(M,(\phi_t)_{t \in \R})$ satisfies the Anosov closing lemma.}

\subsection{Anosov alternative and mixing properties}

{\defn Let $(M,(\phi_t)_{t \in \R})$ be a metric flow. The \textit{length spectrum of the periodic orbits} is the subset $\Lambda \subset \R$ given by the lengths of the periodic orbits of $(M,(\phi_t)_{t \in \R})$: 
$$\Lambda := \lbrace t \in \mathbb{R}_+^\ast \text{ | } \exists a \in M \text{ , } \phi_t(a) = a \rbrace.$$

We say that the length spectrum $\Lambda \subset \R$ of $(M,(\phi_t)_{t \in \R})$ is \textit{arithmetic} if there exists $a \in \R$ such that $\Lambda \subset a.\mathbb{Z}$.}

{\Thm[{\citep[Theorems 3,4 and 6]{Coud}}]\label{AnosovAlternative} Let $(M,(\phi_t)_{t \in \R})$ be a compact metric flow. Suppose that $(M,(\phi_t)_{t \in \R})$ admits a dense set of recurrent points\footnote{Recall that a point $p \in M$ is said to be recurrent (for the flow $(\phi_t)_{t \in \R}$), if for every neighbourhood $U$ of $p$ in $M$ and every  $L > 0$, there exists $t \geq L$ such that $\phi_t(p) \in U$.}, admits a local product structure and satisfies the Anosov closing lemma. The following are equivalent: 
\begin{itemize}
\item[(i)] The flow $(M,(\phi_t)_{t \in \R})$ is weakly topologically mixing.
\item[(ii)] The flow $(M,(\phi_t)_{t \in \R})$ is topologically mixing.
\item[(iii)] The length spectrum of the periodic orbits is not arithmetic.
\end{itemize}}

{\rem Theorem \ref{AnosovAlternative} can be seen as a topological counterpart for the Anosov alternative which states that:
\begin{itemize}
\item The (topological) mixing and (topological) weakly mixing properties are equivalent for Anosov flows. 
\item Any topologically transitive Anosov flow is either mixing or the suspension of a diffeomorphism. In the latter case, the length spectrum of the periodic orbits is arithmetic.
\end{itemize}}

{\rem\label{strategy} Theorem \ref{AnosovAlternative} provides an efficient tool to prove that an Anosov flow $(M,(\phi_t))_{t \in \R}$ is topologically weakly mixing, a key assumption of Theorem B in \cite{moi}. Indeed it shows that, to establish this property, we may proceed as follows:
\begin{itemize}
\item First, prove that the recurrent points are dense in $M$. Notice that this property only depends on the partition of $M$ into orbits of $(\phi_t)_{t \in \R}$. In other words, this property only depends on the foliation by curves defined by the vector field.
\item Then, prove that the length spectrum of the periodic orbits is not arithmetic. In contrast with the preceding case, this is not a property of the partition of $M$ into orbits of  $(\phi_t)_{t \in \R}$  but also depends on the temporal parametrization of the orbits.
\end{itemize}}

{\Thm\label{weaklymixing} Let $(M,g)$ be a compact Riemannian manifold with negative curvature. The geodesic flow $(SM,(g_t)_{t \in \R})$ on the sphere bundle of $M$ is topologically mixing.}

\begin{proof}
By Theorem \ref{anosov}, the unitary geodesic flow $(SM,(g_t)_{t \in \R})$ on the sphere bundle of $M$ is an Anosov flow. The strategy of remark \ref{strategy}, 
reduces the previous theorem to the two following statements: 

{\Thm[{\cite{Hed}}]\label{transitive} Let $(M,g)$ be a compact Riemannian manifold with negative curvature. The recurrent points of the unitary geodesic flow $(SM,(g_t)_{t \in \R})$ on the sphere bundle of $M$ are dense.} 

{\Thm[{\citep[Section II]{Dalbo}}]\label{arithmetic} Let $(M,g)$ be a compact Riemannian manifold with negative curvature. The length spectrum of the periodic orbits of the unitary geodesic flow is not arithmetic.}
\end{proof}

{\rem In the case of a compact Riemannian manifold with constant negative curvature, shorter proofs of Theorem \ref{weaklymixing} can be given by means of representation theory (see for instance \citep[Chapter 4, Section 6]{Glas}). However, in this article, we need to control the topological dynamics of the geodesic flows attached to algebraic metrics, constructed by some approximation procedure (cf. Section 2) and we cannot assume that they have constant negative curvature.} 

\section{Proof of the Main Theorem}

\subsection{Constructing pseudo-Riemannian varieties over $\R$ with negative curvature}
We show the existence of a large supply of smooth pseudo-Riemannian algebraic varieties over $\R$ to which our main theorem (Theorem \ref{Int3}), concerning the orthogonality to the constants of the geodesic flow, will apply. This is actually an easy illustration of the constructions in the second section of this article and of the following two elementary lemmas. 

{\lem\label{irreducible} Let $X$ be a quasi-affine variety over $\R$ and  let $M \subset X(\R)^{an}$ be a connected non-empty $\R$-analytic submanifold.
Then the Zariski-closure $\overline{M}$ of $M$ in $X$ is an irreducible variety over $\R$.} 

\begin{proof}
First, we may suppose that $X$ is affine by working in the Zariski-closure of $X$ relatively to an embedding $i : X \longrightarrow \mathbb{A}^n$. Therefore, we can suppose that the Zariski-closure of $M \subset X(\R)^{an}$ is the Zariski closure of $M \subset \R^n$.

Let $M \subset \R^{an}$ be a connected non-empty analytic submanifold and set $Z = \overline{M}$ . The ideal $I_Z$ of $\R[X_1,\cdots, X_n]$ associated to $Z$ is the kernel of the morphism of rings given by restriction to $M$: 
$$ \R[X_1,\cdots, X_n] \longrightarrow \mathcal O^{an}(M)$$
where $\mathcal O^{an}(M)$ is the ring of analytic functions on $M$. As $M$ is a connected analytic manifold, this ring is an integral domain and its kernel $\mathcal I_Z$ is prime. We conclude that $Z$ is irreducible.
\end{proof}

We will need to work with absolutely irreducible varieties. For that matter, we will  use the following lemma.

{\lem\label{absolutelyirreducible} Let $k$ be a field and $X$ be an irreducible variety over $k$. Suppose that $X(k)$ is Zariski-dense in $X$. Then $X$ is an absolutely irreducible variety over $k$.}

This is well-known. For the sake of completeness, we include a proof.
\begin{proof}
It suffices to prove the lemma for an affine irreducible variety $X$ over $k$ with $X(k)$ Zariski-dense in $X$.

Recall that $X$ is absolutely irreducible if and only if the field extension $k \subset k(X)$ is regular.
As $X(k)$ is Zariski-dense in $X$, there exists an elementary extension $k \preceq k^{\ast}$ and $\overline{c} \in k^n$ realising the generic type of $X$ over $k$.
Every elementary extension of fields is regular and a subextension of a regular extension is also regular.
We conclude that $k \subset k^{\ast}$ is a regular extension and therefore $k \subset k(\overline{c})$ is also regular.

As $\overline{c}$ realises the generic type of $X$, the extensions $k \subset k(\overline{c})$ and $k \subset k(X)$ are isomorphic. In particular, the latter is regular.
\end{proof}

{\Thm\label{montheoremepetit} Let $M$ be a compact $\mathcal{C}^\infty$-manifold which admits a Riemannian metric with negative curvature. There exists an absolutely irreducible smooth pseudo-Riemannian  variety $(X,g)$ over $\R$ such that: 
\begin{itemize}
\item[(i)] The real algebraic set $X(\R)$ is non-singular and Zariski-dense in $X$. Moreover, the smooth manifolds $X(\R)^{an}$ and $M$ are diffeomorphic.
\item[(ii)] The pseudo-Riemannian manifold $(X(\R),g_\R)^{an}$ is a Riemannian manifold with negative curvature. 
\end{itemize}}

\begin{proof}
Let $M$ be a compact manifold. The function $g \mapsto \mathrm{sup}_{P \in \mathcal G_2(M)} K(P)$ which associates to a Riemannian metric $g$, the maximum of the sectional curvature is continuous for the $\mathcal C^2$-topology on $\mathcal C^\infty(M,S^2TM)$.
As $M$ admits a Riemannian metric with negative curvature, by applying theorem \ref{firstapproximation}, there exists a smooth pseudo-Riemannian variety $(X,g)$ over $\R$ such that:
\begin{itemize}
\item[(a)] The real algebraic set $X(\R)$ is non-singular and the smooth manifolds $X(\R)^{an}$ and $M$ are diffeomorphic.
\item[(b)] The pseudo-Riemannian manifold $(X(\R),g_\R)^{an}$ is a Riemannian manifold with negative curvature.
\end{itemize}

Let $Y = \overline{X(\R)} \subset X$ be the Zariski-closure of the real points of $X$ and $\tilde{Y} = Y_{reg}$ the regular points of $Y$. As the real algebraic set $X(\R)$ is non-singular, we have $\tilde{Y}(\R) = X(\R)$ and it is easy to check that $(\tilde{Y}, g_{|\tilde{Y}})$ satisfies the same properties (a) and (b).

Moreover, by definition $\tilde{Y}(\R)$ is Zariski-dense in $Y$ and it only remains to check that $Y$ is an absolutely irreducible variety over $\R$. This is given by lemma \ref{irreducible} together with lemma \ref{absolutelyirreducible}. 
\end{proof}

{\Thm\label{montheoremepetit2} Let $M$ be a compact $\mathcal{C}^\infty$-manifold which admits a Riemannian metric with negative curvature. There exists an absolutely irreducible smooth and quasi-affine variety $X \subset \mathbb{A}^n$ over $\R$ such that:
\begin{itemize}
\item[(i)] The real algebraic set $X(\R)$ is non-singular and $X(\R)^{an}$ is diffeomorphic to $M$.
\item[(ii)] The pseudo-Riemannian manifold $(X(\R),g_\R)^{an}$ is a Riemannian manifold with negative curvature, where $g$ is the restriction of the Euclidean metric on $\mathbb{A}_\R^n$ to $X$. 
\end{itemize}}

The same proof as the proof of Theorem \ref{montheoremepetit} goes through here after applying Theorem \ref{secondapproximation} instead of Theorem \ref{firstapproximation}.

\subsection{The main theorem} 
 
{\Thm\label{montheoreme} Let $(X,g)$ be an absolutely irreducible smooth pseudo-Riemannian variety over $\R$, with Zariski-dense real points, such that $(X(\R),g_\R)^{an}$  is a compact Riemannian manifold.

Then the unitary geodesic $D$-variety associated to $(X,g)$ is absolutely irreducible and, if the sectional curvature of $(X(\R),g_\R)^{an}$ is negative, its generic type is orthogonal to the constants.}

\begin{proof}
Let $(X,g)$ be an absolutely irreducible pseudo-Riemannian variety over $\R$ with 
$X(\R)$ non-singular, Zariski-dense in $X$ and suppose that $(X(\R),g_\R)^{an}$  is a compact and connected Riemannian manifold with negative curvature.

We denote by $(SX,v)$ the associated unitary geodesic $D$-variety and we check that it satisfies the hypothesis of Theorem B in \cite{moi}.
\begin{itemize}
\item It is easy to see that $SX(\R)$ is Zariski-dense in $SX$.

Indeed, recall that $Y(\R)$ is Zariski dense in $Y$ iff $\mathrm{dim}(Y(\R)) = \mathrm{dim}(Y)$ where the first dimension is the semi-algebraic dimension. But we can compute the semi-algebraic dimension of $SX(\R)$ as: 
$$ \mathrm{dim}(SX(\R)) = 2.\mathrm{dim}(X(\R)) - 1 = 2.\mathrm{dim}(X) - 1 = \mathrm{dim}(SX),$$
because $X(\R)$ is Zariski-dense in $X$.

\item Moreover, as $SX(\R)$ is connected and Zariski-dense in $SX$, by applying lemmas \ref{irreducible} and \ref{absolutelyirreducible}, we get that $SX$ is an absolutely irreducible variety over $\R$.

\item We now show that the real flow of the unitary geodesic $D$-variety $(SX,v)$ is weakly mixing.
From Lemma \ref{realflow}, we know that this flow is the unitary geodesic flow of the  compact Riemannian manifold $(X(\R),g_\R)^{an}$. By hypothesis, its sectional curvature is negative. Consequently, using Theorem \ref{weaklymixing}, we conclude that the unitary geodesic flow is mixing and therefore weakly mixing. 
\end{itemize}
By applying Theorem B in \cite{moi} to the compact set $SX(\R)$, we get that the generic type of the absolutely irreducible $D$-variety $(SX,v)$ is orthogonal to the constants.
\end{proof}

{\rem In the proof of Theorem \ref{montheoreme}, we have applied Theorem B of \cite{moi} to all real points of $X$, whereas it suffices to know the dynamics on a compact Zariski-dense subset $K \subset X(\R)$ to conclude that the generic type is orthogonal to the constants.

In particular, we can get a stronger form of the theorem above by only requiring that a connected and compact component $C \subset X(\R)$  is non-singular and that the pseudo-Riemannian manifold $(C,g_C^{an})$ is a Riemannian manifold with compact curvature.

With this stronger version at hand, one can bypass the use of the Tognoli's Theorem (Theorem \ref{Tognoli}) in the existence statements \ref{montheoremepetit} and \ref{montheoremepetit2} and replace it by its weaker version (see Remark \ref{Nashtheorem}).}

{\cor Let $(X,g)$ be an absolutely irreducible smooth pseudo-Riemannian variety over $\R$, with Zariski-dense real points, such that $(X(\R),g_\R)^{an}$  is a compact Riemannian manifold.

If the sectional curvature of this Riemannian manifold is negative, then every non zero fibre of the morphism of $D$-varieties
$$ H_g : (T_{X/k} ,X_g) \longrightarrow (\mathbb{A}^1,0)$$
is an absolutely irreducible $D$-variety with  generic type orthogonal to the constants.}

It suffices to apply Theorem \ref{montheoreme} and then Proposition \ref{fiberchange}. Note that by Theorem \ref{montheoremepetit}, there exist many examples of such pseudo-Riemannian varieties.

{{\rem Let $(X,g)$ be an absolutely irreducible smooth pseudo-Riemannian variety over $\R$, with Zariski-dense real points, such that $(X(\R),g_\R)^{an}$  is a compact and connected Riemannian manifold.
The morphism of $D$-varieties $$ H_g : (T_{X/k} ,X_g) \longrightarrow (\mathbb{A}^1,0)$$ has three kinds of non-zero fibres: 
\begin{itemize}
\item The fibres over positive real numbers which are real $D$-varieties with Zariki-dense real points.
\item The fibres over negative real numbers which are real $D$-varieties without any real points.
\item The fibres over complex non real numbers which are complex $D$-varieties.
\end{itemize}
In particular, thanks to Proposition \ref{fiberchange}, we can obtain orthogonality results for some real $D$-varieties without real points and even for some complex (non real) $D$-varieties.}

{\cor\label{theoremA} Let $M \subset \R^n$ be a non-empty compact connected and non-singular algebraic subset of the Euclidean space $\R^n$.

The system $(S)$ of differential equations describing the motion of a particle in the Euclidean space $\R^n$ with non-zero fixed energy, constrained to move without friction along the submanifold $M$ is an absolutely irreducible system of algebraic differential equations. 

If the restriction of the Euclidean metric to $M$ has negative curvature, then the generic type of the system $(S)$ is orthogonal to the constants.}

\begin{proof}
Set $M = X(\R)$ with $X$ a quasi-affine variety over $\R$. We can suppose (after restricting to an open subset of $X$) that $X$ is non singular and that the restriction to $X$ of the standard Euclidean metric $g$ on $\mathbb{A}^n_\R$ is non-degenerate.

The system $(S)$ with fixed energy $E_0 \in \mathbb{C}$ is the system of differential equations given as the fibre (as a $D$-variety) over $E_0$ of the morphism of $D$-varieties 
$$ H_g : (T_{X/k} ,X_g) \longrightarrow (\mathbb{A}^1,0).$$ It suffices to apply Theorem \ref{montheoreme} to $(X,g)$ and then Proposition \ref{fiberchange} to conclude that any non-zero fibre is absolutely irreducible and that its generic type is orthogonal to the constants. 
\end{proof}

\bibliographystyle{alpha}
\bibliography{bibliographie}

\end{document}